\documentclass{article}
\usepackage{amssymb,amsmath,amsthm,graphicx}

\def\tr{\triangleright}

\newtheorem{theorem}{Theorem}

\theoremstyle{definition}
\newtheorem{example}{Example}
\newtheorem{definition}{Definition}
\newtheorem{remark}{Remark}

\date{}

\title{\Large \textbf{Quandle Cocycle Quivers}}

\author{Karina Cho\footnote{Email: kcho@hmc.edu.}
\and
Sam Nelson\footnote{Email: Sam.Nelson@cmc.edu. Partially supported by Simons Foundation collaboration grant 316709}}

\begin{document}
\maketitle

\begin{abstract} We incorporate quandle cocycle information into the
quandle coloring quivers we defined in \cite{CN1} to define weighted
directed graph-valued invariants of oriented links we call \textit{quandle
cocycle quivers}. This construction turns the quandle cocycle invariant
into a small category, yielding a categorification of the quandle cocycle
invariant. From these graphs we define several new link invariants including
a 2-variable polynomial which specializes to the usual quandle cocycle 
invariant. Examples and computations are provided.
\end{abstract}

\parbox{5.25in}{\textsc{Keywords:} Quandles, Enhancements, Quivers,
Cocycle Enhancements

\smallskip

}

\section{\large\textbf{Introduction}}\label{I}

\textit{Quandles}, algebraic structures whose axioms are derived from the
Reidemeister moves, were introduced in \cite{J,M} and have been studied
ever since. Associated to every oriented knot or link there exists a  
\textit{fundamental quandle}, also called the \textit{knot quandle} of
the knot or link. The isomorphism class of the fundamental quandle of a knot
determines the group system and hence the knot up to mirror image; in this
sense the knot quandle is a complete invariant for knots, though not for
split links or other generalizations such as virtual knots.
The number of homomorphisms from a knot's fundamental quandle
to a finite coloring quandle, called the \textit{quandle counting invariant},
is an integer-valued invariant of oriented knots and links.

The quandle counting invariant
forms the basis for a class of computable knot and link invariants known
as \textit{enhancements}, including in particular \textit{quandle cocycle 
invariants}. In \cite{CJKLS}, a theory of quandle cohomology was introduced.
Quandle-colored link diagrams can be identified with certain elements of
the second homology of the coloring quandle, and it naturally follows that
elements of the second cohomology define invariants of quandle-colored
isotopy. See \cite{CJKLS,EN} etc. for more. 

In \cite{CN1}, we introduced \textit{quandle coloring quivers}, directed 
graph-value invariants of oriented knots and links associated to
pairs $(X,S)$ consisting of a finite quandle $X$ and a subset $S$ of the 
ring of endomorphisms of $X$. From these quivers we defined a new polynomial
knot invariant, using this quiver structure to enhance the quandle counting
invariant, called the \textit{in-degree polynomial}. Since directed graphs
can be interpreted as categories, this construction yielded a categorification
of the quandle counting invariant.

In this paper we enhance quandle coloring quivers with quandle cocycles, 
obtaining \textit{quandle cocycle quivers.} As in the previous case, this
construction yields a categorification of the quandle cocycle invariant. From
these quivers we derive new polynomial knot and link invariants, including
a two-variable enhancement of the quandle 2-cocycle invariant which 
has both the quandle counting invariant and the quandle cocycle invariant
as specializations but in general is a stronger invariant than either. The 
paper is organized as follows. In Section \ref{QC} we review quandles and 
quandle cohomology. In Section \ref{QQ} we define quandle cocycle quivers and
provide examples demonstrating the computation of the invariant and showing
that the new polynomial invariant is a proper enhancement. We conclude in 
Section \ref{Q} with some questions for future research.

\section{\large\textbf{Quandles and Quandle Cohomology}}\label{QC}
We begin with some preliminary definitions.
%start quandle basics
\begin{definition}
A set $X$ with a binary operation $\tr$ is a \textit{quandle} if it satisfies
\begin{enumerate}
    \item[(i)] for all $x\in X$, $x\tr x=x$ 
    \item[(ii)] for all $y\in X$, the map $f_y: X\to X$ defined by $f_y(x)=x\tr y$ is a bijection, and
    \item[(iii)] for all $x, y, z\in X$, $(x\tr y)\tr z = (x\tr z)\tr(y\tr z)$ .
\end{enumerate}
\end{definition}

\begin{definition}
For an oriented link $L$, associate a label to each arc of $L$. At each crossing of $L$, assign a relation between arc labels as in the figure below.
\[\includegraphics[scale=0.3]{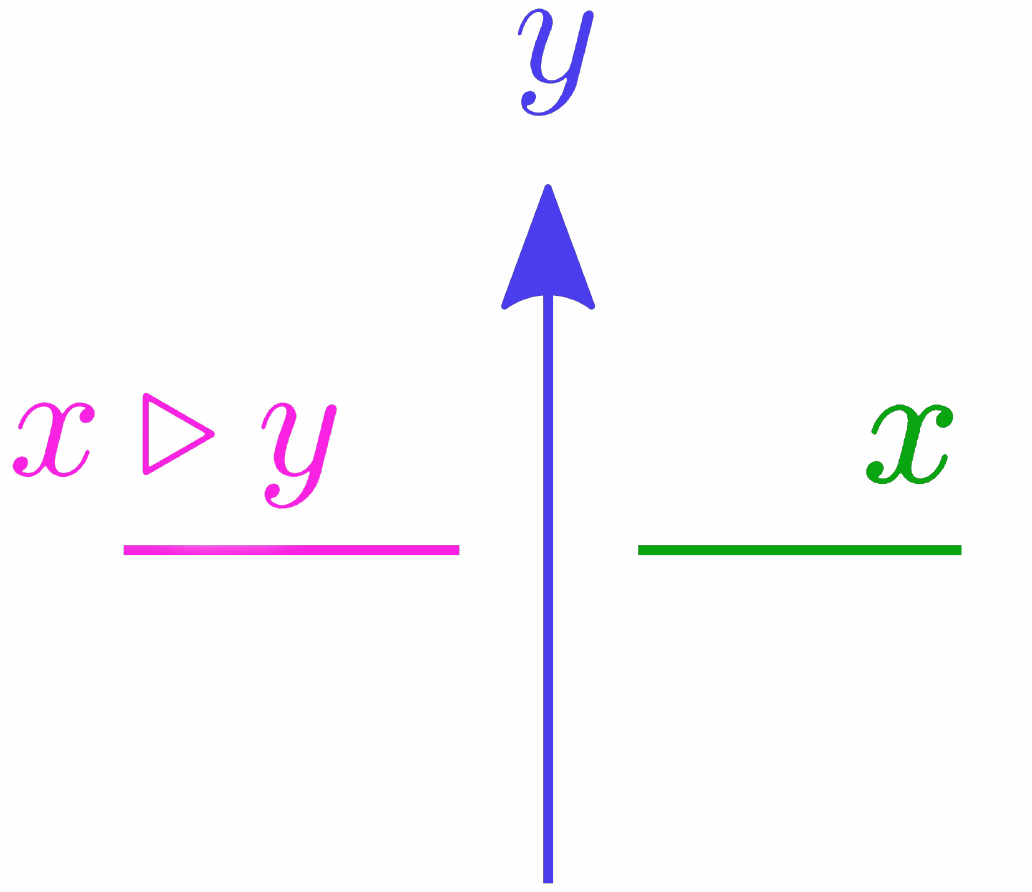}
\]
In words, if we orient a crossing so that the overstrand (labeled $y$) points up, and we have arcs labeled $z$ and $x$ on the left and right respectively of the arc labeled $y$, then we have the crossing relation $z=x\tr y$.

The \textit{fundamental quandle} $Q(L)$ is the quandle generated by the set of arc labels under the equivalence relations defined by crossing relations for each crossing of $L$.
\end{definition}

\begin{remark}
The quandle axioms are defined in such a way that the fundamental quandle of an oriented link is invariant under Reidemeister moves, so $Q(L)$ is an invariant of $L$.
\end{remark}

\begin{definition}
Let $X$ by a finite quandle and $L$ be an oriented link. Then a homomorphism $\phi: Q(L)\to X$ is called an \textit{$X$-coloring} of $L$.
\end{definition}

\begin{example}
Let $X$ be the quandle defined by the operation table
\[\begin{array}{r|rrrrrr}
\tr & 1 & 2 & 3 & 4  \\\hline
    1 & 1 & 3 & 4 & 2 \\
    2 & 4 & 2 & 1 & 3 \\
    3 & 2 & 4 & 3 & 1 \\
    4 & 3 & 1 & 2 & 4 
\end{array}\]
and $L$ be the figure-8 knot as drawn below, with arc labels $a,b,c,d$. Then 
\[Q(L)=\langle a,b,c,d\ |\ b=a\tr d=d\tr c, c=a\tr b=d\tr a \rangle,\] 
and $\phi: Q(L)\to X$ that maps $a\mapsto 1, b\mapsto 2, c\mapsto 3, d\mapsto 4$ is an $X$-coloring of $L$, which we may visualize with the picture below.

\[%\includegraphics[scale=0.25]{coloringfigure8.pdf}
\includegraphics[scale=0.25]{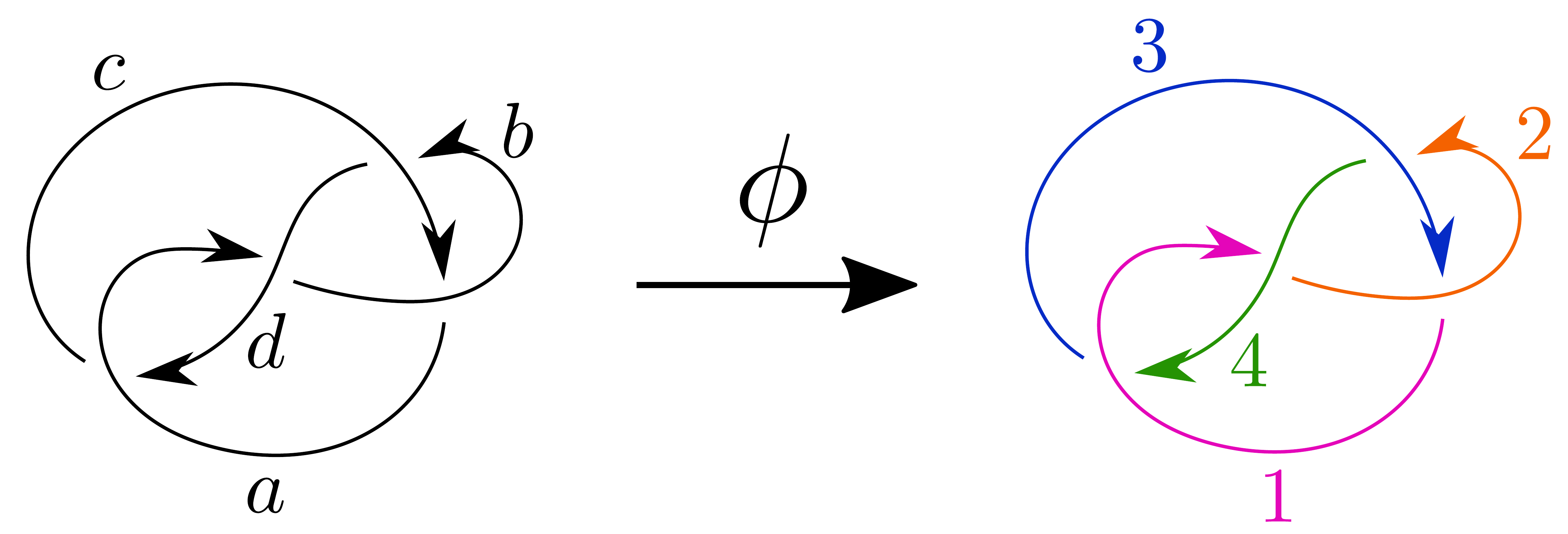}
\]

\end{example}
%end quandle basics

The quandle axioms are the conditions required to ensure that for a given
quandle coloring of a diagram on one side of a move, there is a unique 
quandle coloring of the diagram on the other side of the move which agrees with the original coloring outside the neighborhood of the move. From this, we 
obtain the following standard result:

\begin{theorem}
Let $X$ be a finite quandle. The number of colorings of an oriented knot
or link diagram is an integer-valued invariant of ambient isotopy.
\end{theorem}

This theorem also follows from the observation that the set of quandle 
colorings of a knot or link diagram can be identified with the set
$\mathrm{Hom}(\mathcal{Q}(L),X)$ of quandle homomorphisms from the 
fundamental quandle of the knot or link $L$ to the coloring quandle $X$.
See \cite{EN} for more.

\begin{definition}
Let $X$ be a finite quandle, $A$ an abelian group and 
$C_n^R(X;A)=A[X^n]$, the set of $A$-linear combinations of ordered $n$-tuples
of $X$. Let $C_n^D(X;A)$ be the subgroup generated by elements 
$(x_1,\dots, x_n)$ with $x_j=x_{j+1}$ for some $j$, and let 
$C_n^Q(X;A)=C_n^R(X;A)/C_n^D(X;A)$. 
Define $\partial_n:C_n^R(X;A)\to C_{n-1}^R(X;A)$ by setting
\[\partial_n(\vec{x})
=\sum_{k=1}^n (-1)^k\left(\partial_n^0(\vec{x})-\partial_n^1(\vec{x})\right)\]
where we have
\begin{eqnarray*}
\partial_n^0(x_1,\dots, x_n) & = & (x_1,\dots, x_{k-1}, x_{k+1},\dots, x_n) \\
\partial_n^0(x_1,\dots, x_n) & = & (x_1\tr x_k,\dots, x_{k-1}\tr x_k, x_{k+1},\dots, x_n)
\end{eqnarray*}
and extending linearly.
Since $\partial_n(C_n^D(X;A))\subset C^D_{n-1}(X;A)$, $\partial_n$ induces 
$\partial_n^Q:C_n^Q(X;A)\to C_{n-1}^Q(X;A)$.
Then the \textit{$n$th quandle homology} of $X$ is
\[H_n^{Q}(X)=\mathrm{Ker} (\partial^Q_n)/\mathrm{Im} (\partial_{n+1}^Q).\]
Dualizing, we have $C_*^n(X;A)=\mathrm{Hom}(C_n^*(X;A),A)$,
$\delta^n:C_R^n(X)\to C_R^{n+1}(X;A)$
defined by 
\[(\delta^n f)(x_1,\dots, x_n)
= f\partial_{n+1}(x_1,\dots, x_{n+1})
\]
and \textit{$n$th quandle cohomology} of $X$ given by
\[H_Q^n(X)=\mathrm{Ker} (\delta_Q^n)/\mathrm{Im} (\delta_Q^{n-1}).\]
%adding 1/31
An element of $\mathrm{Ker} (\delta_Q^n)$ is called a 
\textit{quandle $n$-cocycle} and 
an element of $\mathrm{Im} (\delta_Q^{n-1})$ is called a 
\textit{quandle $n$-coboundary}. %check this
\end{definition}

%added by karina 1/31
In this paper, we are particularly interested in quandle 2-cocycles, which are 
maps $\phi:A[X \times X] \to A$ for an abelian group $A$ 
(usually $\mathbb{Z}_n$ for us). These can be written as linear 
combinations of elementary functions $\chi_{i,j}:X \times X \to A$ where
\[
    \chi_{i,j}(x_1,x_2) = 
    \begin{cases}
        1, & \text{for } i=x_1, j=x_2\\
        0, & \text{otherwise.} 
    \end{cases}  
\]
\begin{remark}\label{rem:BWcond}
A quick computation will show that $\phi$ is a quandle 2-cocycle if and only if $\phi$ satisfies the condition that 
\[\phi(x,y)+\phi(x\tr y,z) = \phi(x,z)+\phi(x\tr z,y\tr z)\] 
for all $x,y,z \in X$ and 
$\phi(x,x) = 0$ for all $x\in X$. %from modding by the degeneracy module
\end{remark}

\begin{definition}
Let $L$ be an oriented link, $X$ be a finite quandle, and $\phi$ a quandle 2-cocycle. Let $v$ be an $X$-coloring of $L$. For an $X$-colored crossing $c$ of $L$, we define $\phi(c)$ by the following based on the whether the crossing is positively or negatively oriented:
\[%\includegraphics[scale=0.2]{boltzmannCrossing2.pdf}
\includegraphics[scale=0.2]{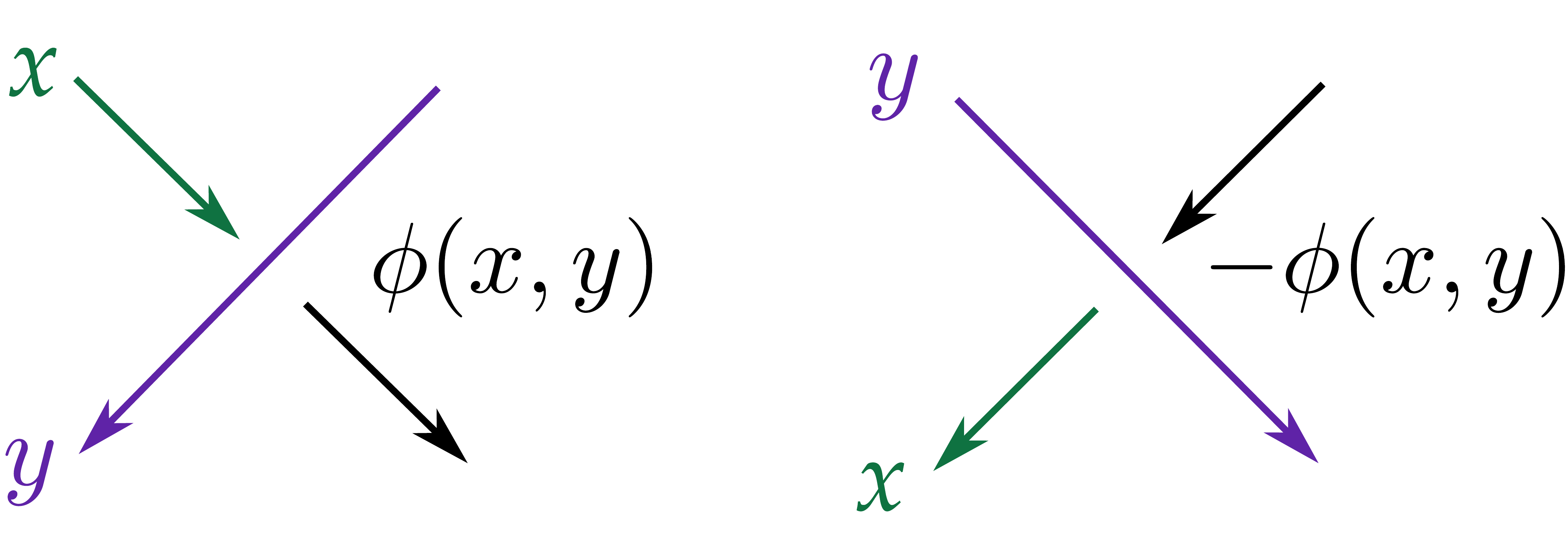}
\]
In the illustration above, $x$ and $y$ are the elements of $X$ determined by 
the coloring $v$. Then we define $\phi(v)$ to be
\[ \phi(v) = \sum_{c \in C}\phi(c), \]
where $C$ is the set of crossings in $L$.
\end{definition}

\begin{remark}
It can be shown that the evaluation of $\phi$ on a quandle colored link 
satisfying the condition described in Remark \ref{rem:BWcond} is equivalent 
to $\phi$ being invariant under a quandle-colored Reidemeister type III move. 
In fact, $\phi$ is an invariant of quandle-colored links. See \cite{CJKLS,EN} 
for more.
\end{remark}

\begin{definition}
Let $X$ be a finite quandle, $L$ an oriented link and 
$\phi\in H^2_Q(X;A)$. Then the polynomial
\[\Phi_X^{\phi}(L)=\sum_{v\in\mathcal{C}(L,X)} s^{\phi(v)}\]
where $\mathcal{C}(L,X)$ is the set of $X$-colorings of a diagram of $L$
and $\phi(v)$ is the Boltzmann weight of the $X$-colored diagram $v$
is the \textit{quandle 2-cocycle invariant} of $L$. See \cite{CJKLS,EN} for
more.
\end{definition}

\section{\large\textbf{Coloring Quivers and Cocycle Quivers}}\label{QQ}

Let $L$ be an oriented link diagram.
In \cite{CN1} we defined the \textit{quandle coloring quiver} invariant
$\mathcal{Q}_X^S(L)$ in the following way: given a finite quandle $X$
and set $S\subset \mathrm{Hom}(X,X)$ of quandle endomorphisms, we make
a directed graph with a vertex for each $X$-coloring of $L$ and a directed
edge from $v_j$ to $v_k$ whenever $v_k=f(v_j)$ in the sense that each 
arc color in $v_k$ is obtained from the corresponding arc color in $v_j$
by applying $f$ for some $f\in S$.

\begin{example}\label{ex1}
The links $L7n1$ and $L7n2$ have isomorphic quandle coloring quivers with
respect to the quandle and endomorphism
\[\includegraphics[scale=0.2]{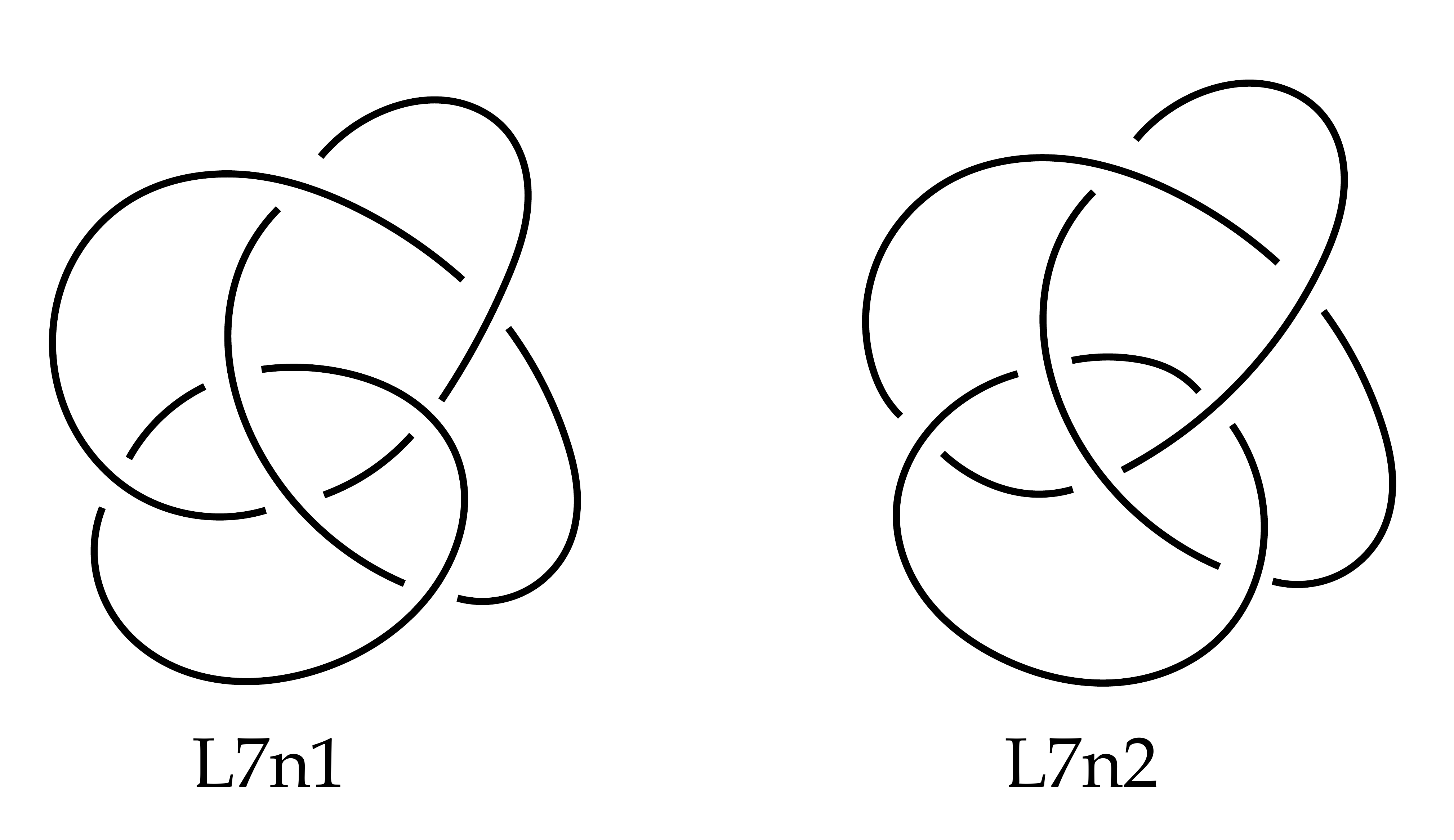}\]
\[
\begin{array}{r|rrrr}
\tr & 1 & 2 & 3 & 4 \\ \hline
1 & 1 & 1 & 1 & 1 \\
2 & 4 & 2 & 2 & 2 \\
3 & 3 & 3 & 3 & 3 \\
4 & 2 & 4 & 4 & 4
\end{array}\quad\quad\quad
f(1)=4, f(2)=f(3)=f(4)=3
\]
namely
\[\includegraphics[scale=1]{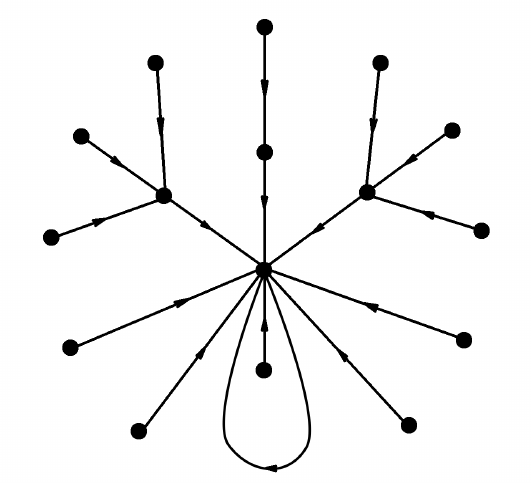}.
\]
\end{example}

\begin{definition}
Let $L$ be an oriented link, $X$ be a finite quandle, $S$ a set of 
quandle endomorphisms from $X$ to $X$, and $\phi$ a 2-cocycle in $C^2_Q(X)$.
Then the \textit{quandle cocycle quiver} $\mathcal{Q}_X^{S,\phi}(L)$ 
is the directed graph with 
vertices corresponding to $X$-colorings of $L$, edges from $v_j$ to
$v_k$ whenever $v_k=f(v_j)$ for some $f\in S$, and weights $\phi(v_j)$ at
each vertex. When $S=\{f\}$ is a singleton we will write $f$ instead of 
$\{f\}$ for simplicity.
\end{definition}

\begin{example}\label{ex2}
The links in example \ref{ex1} are not distinguished by their coloring
quivers with respect to the given quandle and endomorphism; however, the
quandle cocycle quiver with cocycle 
\[\phi=\chi_{1,2}+2\chi_{1,3}+\chi_{1,4}+2\chi_{2,1}+3\chi_{3,2}+3\chi_{3,4}+\chi_{4,1}\in C_Q^2(X;\mathbb{Z}_4)\]
does distinguish the links.
\[\begin{array}{cc}
\includegraphics{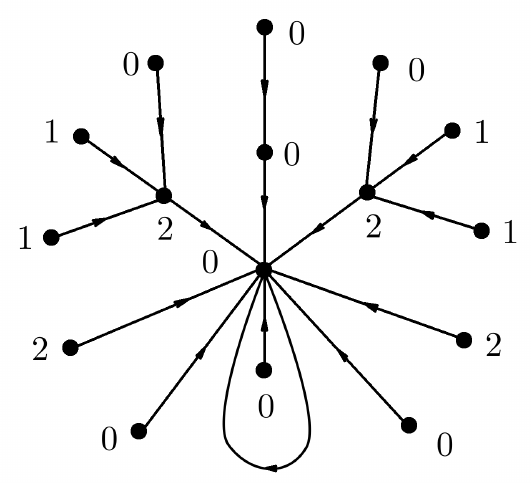} & \includegraphics{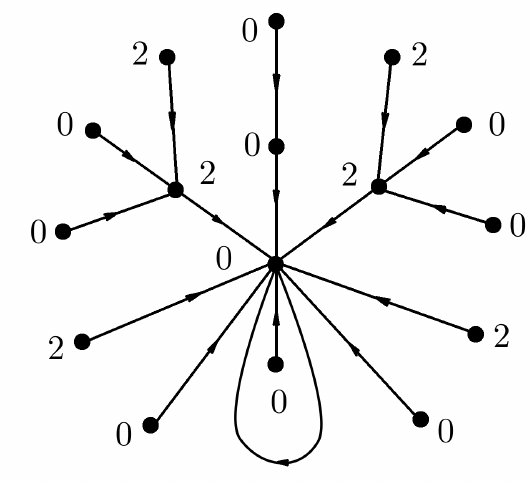}\\
\mathcal{Q}_X^{f,\phi}(L7n1) &  \mathcal{Q}_X^{f,\phi}(L7n2)
\end{array}\]
\end{example}

\begin{definition}
Let $L$ be a link, $X$ a finite quandle, $S\subset \mathrm{Hom}(X,X)$ and
$\phi\in C^2_Q(X;A)$. We define the \textit{quiver enhanced cocycle polynomial}
to be the polynomial
\[\Phi_X^{S,\phi}(L)=\sum_{e\in E(\mathcal{Q}_X^S(L))} s^{\phi(v_j)}t^{\phi(v_k)}\]
where the edge $e$ is directed from vertex $v_j$ to vertex $v_k$ in the quandle
coloring quiver $\mathcal{Q}_X^S(L)$. 
\end{definition}

\begin{remark}
Since each directed edge contributes its $s^{\phi(v_j)}t^{\phi(v_k)}$ value to the
polynomial  $\Phi_X^{S,\phi}(L)$ independently, if we regard the quandle coloring
quiver $\mathcal{Q}_{X,S}(L)$ as the union of the quivers 
$\mathcal{Q}_{X,f}(L)$ for endomorphisms $f\in S$, then the polynomial
can be separated into a sum of the cocycle quiver polynomials for each
individual endomorphism:
\[\Phi_X^{S,\phi}(L)=\sum_{f\in S} \Phi_X^{f,\phi}(L).\]
It follows that
evaluating $\Phi_X^{S,\phi}(L)$ at $t=1$ yields $|S|\Phi_X^{\phi}(L)$. In 
particular, when $|S|=1$, $\Phi_X^{S,\phi}(L)$ evaluates at $t=1$ to the 
classical quandle 2-cocycle invariant as defined in \cite{CJKLS} 
(see e.g. Example 8 in \cite{CJKLS} and note that our $s$ is their $t$).

Similarly, if $f$ is the identity endomorphism, then 
$\Phi_X^{f,\phi}(L)$ is the quandle cocycle invariant evaluated at $st$.
\end{remark}

\begin{example} \label{ex6}
In \cite{CN1} Example 6 we gave an example of two links $L6a1$ and $L6a5$
which have the same counting invariant value $\Phi_X^{\mathbb{Z}}(L)=16$ with 
respect to the quandle $X$ with operation table
\[\begin{array}{r|rrrr}
\tr & 1 & 2 & 3 & 4 \\ \hline
1 & 1 & 3 & 1 & 3 \\
2 & 4 & 2 & 4 & 2 \\
3 & 3 & 1 & 3 & 1\\
4 & 2 & 4 & 2 & 4
\end{array}
\]
but are distinguished by the isomorphism class of their quandle coloring
quivers $\mathcal{Q}_X^{f}(L)$ where $f:X\to X$ is the endomorphism
given by $f(1)=f(3)=4$ and $f(2)=f(4)=2$. (The published version incorrectly
lists this as ``$f(1)=4, f(1)=2, f(1)=4, f(1)=2$''.) We note that for any 
coboundary $\phi$, the Boltzmann weights are all 0 and the counting invariant
$\Phi_X^{\mathbb{Z}}(L)=16$
is equal to the quiver enhanced polynomial, so this example also shows that the
Boltzmann weight enhanced quiver is not determined by the quiver enhanced 
polynomial.
\end{example}

The links in example \ref{ex2} are distinguished by their quandle cocycle 
quivers, but they are already distinguished by their quandle cocycle invariants
with respect to the given quandle and cocycle. The next example shows that the 
quandle cocycle quiver can distinguish knots which have the same
quandle cocycle invariant.

\begin{example} \label{ex3}
Let $X$ be the quandle defined by operation table 
\[\begin{array}{r|rrrrrr}
\tr & 1 & 2 & 3 & 4 & 5 & 6 \\\hline
1 & 1 & 3 & 2 & 5 & 4 & 1 \\
2 & 3 & 2 & 1 & 6 & 2 & 4 \\
3 & 2 & 1 & 3 & 3 & 6 & 5 \\
4 & 5 & 6 & 4 & 4 & 1 & 2 \\
5 & 4 & 5 & 6 & 1 & 5 & 3 \\
6 & 6 & 4 & 5 & 2 & 3 & 6.
\end{array}\]
This quandle has endomorphism $f(1)=f(6)=2$, $f(2)=f(5)=4$ and $f(3)=f(4)=6$
and cocycle 
\begin{eqnarray*}
\phi & = & 
2\chi_{(1,2)}+
2\chi_{(1,3)}+
2\chi_{(1,4)}+
2\chi_{(1,5)}+
\chi_{(2,3)}+
2\chi_{(2,4)}+
\chi_{(3,2)}+
2\chi_{(3,5)}\\
& & +
2\chi_{(4,2)}+
\chi_{(4,5)}+
2\chi_{(5,3)}+
\chi_{(5,4)}+
\chi_{(6,2)}+
\chi_{(6,3)}+
\chi_{(6,4)}+
\chi_{(6,5)}
\end{eqnarray*}
in $C^2_Q(X;\mathbb{Z}_3)$. Then the knots $6_1$ and $7_7$ both have quandle
cocycle polynomial
\[\Phi_X^{\phi}(6_1)=6+12s+12s^2=\Phi_X^{\phi}(7_7)\]
but are distinguished by their quiver enhanced polynomials
\[\Phi_X^{f,\phi}(6_1)=6+12st+12s^2t^2 \ne
6+12st^2+12s^2t=\Phi_X^{f,\phi}(7_7). \]
\end{example}

\begin{example}
Continuing with the same quandle from example \ref{ex3},
we computed $\Phi_X^{f,\phi}(L)$ using our \texttt{python} code for the 
prime knots with up to eight crossings and prime links with up to seven 
crossings with the  cocycle with $\mathbb{Z}_4$ coefficients
\[\begin{array}{rcl} \phi & = & 
\chi_{13}+3\chi_{14}+2\chi_{15}+3\chi_{21}+3\chi_{23}+2\chi_{24}+\chi_{25}\\
& & +\chi_{31}+3\chi_{35}+3\chi_{36}+\chi_{41}+\chi_{42}+2\chi_{45}+3\chi_{46} \\
& & +3\chi_{51}+\chi_{54}+\chi_{56}+3\chi_{62}+3\chi_{64}+\chi_{65}
\end{array}\]
and three arbitrarily chosen
endomorphisms (where we write $f$ by specifying $[f(1),\dots, f(n)]$):
\[\begin{array}{rcl}
f_1 & = & [1,2,3,3,2,1]\\
f_2 & = & [1,5,4,3,2,6]\\
f_3 & = & [3,1,2,5,6,4]\\
\end{array}.\] The results are collected in the table. For simplicity we list 
only the nontrivial values; unlisted knots have 
$\Phi_X^{f,\phi}(L)=6$.
\[
\begin{array}{c|ccc}
L & \Phi_X^{f_1,\phi}(L) & \Phi_X^{f_2,\phi}(L) & \Phi_X^{f_3,\phi}(L) \\ \hline
3_1 & 30 & 30 & 30 \\
6_1 & 12s^2t^2+4s^2+4t^2+10 & 16s^2t^2+14 & 10s^2t^2+6s^2+6t^2+8\\
7_4 & 30 & 30 & 30 \\
7_7 & 12s^2t^2+4s^2+4t^2+10 & 16s^2t^2+14 & 10s^2t^2+6s^2+6t^2+8\\
8_{10} & 54 & 54 & 54 \\
8_{11} & 12s^2t^2+4s^2+4t^2+10 & 16s^2t^2+14 & 10s^2t^2+6s^2+6t^2+8\\
8_{15} & 54 &  54 & 54 \\
8_{18} & 48s^2t^2+16s^2+16t^2+22 & 64s^2t^2+38 & 40s^2t^2+24s^2+24t^2+14 \\
8_{19} & 54 & 54 &54 \\
8_{20} & 54 & 54 & 54 \\
8_{21} & 54 & 54 & 54\\
L2a1 & 12 & 12 & 12\\
L4a1 & 12 & 12 & 12\\
L5a1 & 12 & 12 & 12 \\
L6a1 & 12s^2t^2+4s^2+4t^2+16 & 16s^2t^2+20 & 10s^2t^2+6s^2+6t^2+14 \\
L6a2 & 12 & 12 & 12\\
L6a3 & 36 & 36 & 36 \\
L6a4 & 24 & 24 & 24 \\
L6a5 & 48 & 48 & 48 \\
L6n1 & 24 & 24 & 24 \\ 
L7a1 & 16s^2t^2+8s^2+16t^2+20 & 24s^2t^2+36 & 12s^2t^2+12s^2+12t^2+24 \\
L7a2 & 12 & 12 &12 \\
L7a3 & 12 & 12 & 12 \\
L7a4 & 12 & 12 & 12 \\
L7a5 & 12s^2t^2+4s^2+4t^2+16 & 16s^2t^2+20 & 10s^2t^2+6s^2+6t^2+14\\
L7a6 & 12 & 12 & 12 \\
L7a7 & 24 & 24 & 24 \\
L7n1 & 12 & 12 & 12 \\
L7n2 & 12 & 12 & 12 \\
\end{array}
\]

\end{example}

Many other link invariants can be defined from quandle cocycle quivers.
For instance, we can modify the incidence matrix of the graph incorporating
the cocycle information:

\begin{definition}
Let $X$ be a quandle, $\phi\in H^{2}_Q(X;A)$ a quandle 2-cocycle with values
in an abelian group $A$, and $S\subset\mathrm{Hom}(X,X)$ a set of endomorphisms.
Then for a oriented link $L$ and choice of numbering for vertices and edges in 
$\mathcal{Q}_X^{f}(L)$, we define a matrix $M_{\mathcal{Q}_X^{f}(L)}$ whose entry in
row $j$ column $k$ is 
\[\left\{\begin{array}{rl}
-\phi(v_j) & v_j=\mathrm{Source}(e_k)\\
\phi(v_j) & v_j=\mathrm{Target}(e_k)\\
0 & \mathrm{Else}
\end{array}\right.
\] 
\end{definition}

The matrix $M_{\mathcal{Q}_X^{f}(L)}$ itself depends on our choice of numbering
for vertices and edges, but we can obtain from it several link invariants
including but not limited to:
\begin{itemize}
\item The rank of $M_{\mathcal{Q}_X^{f}(L)}$,
\item The isomorphism class of the linear transformation between $R$-modules
determined by the matrix,
\item The Smith normal form of the matrix when $R$ is a PID,
\item The eigenvalues and characteristic polynomial when $M_{\mathcal{Q}_X^{f}(L)}$
is a square matrix,
\item The elementary ideals of $M_{\mathcal{Q}_X^{f}(L)}$
\end{itemize}
and more.

\begin{example}\label{ex4}
Let $L$ be the link $L4a1$ (the $(4,2)$ torus link) and consider the quandle,
endomorphism and cocycle
\[
\begin{array}{r|rrr}
\tr & 1 &2 & 3\\ \hline
1 & 1 & 1 & 2 \\
2 & 2 & 2 & 1 \\
3 & 3 & 3 & 3
\end{array},\quad 
f=[1,1,2], \quad 
\phi=2\chi_{13}+3\chi_{23}+4\chi_{31}+4\chi_{32}\in H^2_Q(X;\mathbb{Z}_5).
\]
Then we compute the matrix $M_{\mathcal{Q}_X^{f}(L)}\in M_9(\mathbb{Z}_5)$ for 
the cocycle quiver:
\[\raisebox{-1in}{\includegraphics{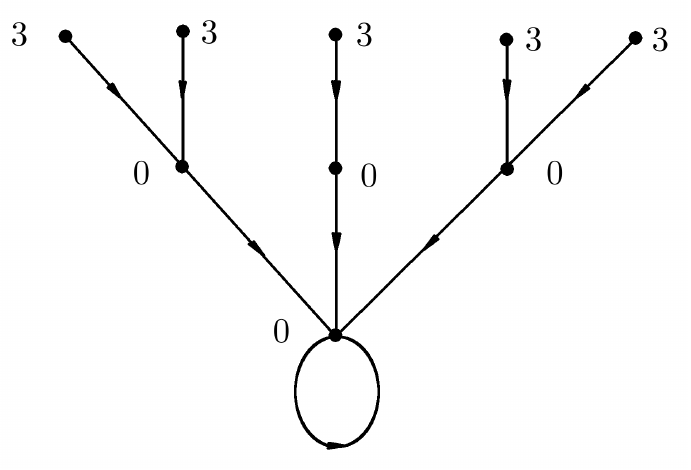}} \quad 
\left[\begin{array}{rrrrrrrrr}
2 & 0 & 0 & 0 & 0 & 0 & 0 & 0 & 0 \\
0 & 2 & 0 & 0 & 0 & 0 & 0 & 0 & 0 \\
0 & 0 & 2 & 0 & 0 & 0 & 0 & 0 & 0 \\
0 & 0 & 0 & 2 & 0 & 0 & 0 & 0 & 0 \\
0 & 0 & 0 & 0 & 2 & 0 & 0 & 0 & 0 \\
0 & 0 & 0 & 0 & 0 & 0 & 0 & 0 & 0 \\
0 & 0 & 0 & 0 & 0 & 0 & 0 & 0 & 0 \\
0 & 0 & 0 & 0 & 0 & 0 & 0 & 0 & 0 \\
0 & 0 & 0 & 0 & 0 & 0 & 0 & 0 & 0 \\
\end{array}\right].\]
Then for example, we obtain cocycle quiver characteristic polynomial
value $(x+3)^5x^4\in \mathbb{Z}_5[x]$.
\end{example}

\begin{example}
Using the same quandle and cocycle as in example \ref{ex4}
with endomorphism $f=[2,1,3]$, the 
link $L7a3$ has cocycle quiver and matrix
\[\raisebox{-1in}{\includegraphics{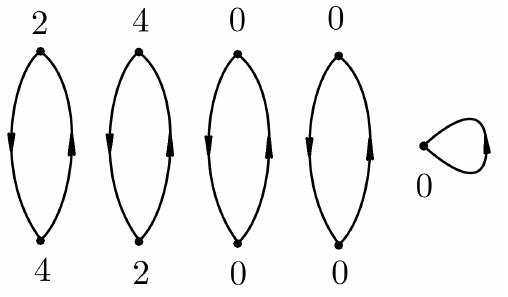}} \quad 
\left[\begin{array}{rrrrrrrrr}
1 & 3 & 0 & 0 & 0 & 0 & 0 & 0 & 0 \\
2 & 4 & 0 & 0 & 0 & 0 & 0 & 0 & 0 \\
0 & 0 & 1 & 3 & 0 & 0 & 0 & 0 & 0 \\
0 & 0 & 2 & 4 & 0 & 0 & 0 & 0 & 0 \\
0 & 0 & 0 & 0 & 0 & 0 & 0 & 0 & 0 \\
0 & 0 & 0 & 0 & 0 & 0 & 0 & 0 & 0 \\
0 & 0 & 0 & 0 & 0 & 0 & 0 & 0 & 0 \\
0 & 0 & 0 & 0 & 0 & 0 & 0 & 0 & 0 \\
0 & 0 & 0 & 0 & 0 & 0 & 0 & 0 & 0 \\
\end{array}\right]\]
with characteristic polynomial $(x^2+3)^2x^5\in \mathbb{Z}_5[x]$.
\end{example}

\section{\large\textbf{Questions}}\label{Q}

We conclude with a few questions and direction for future work.

As we saw in example \ref{ex6} the cocycle polynomial does 
not determine the cocycle quiver. We are curious about which $R$-colored
quivers are obtainable as quandle cocycle quivers of knots and links.
For example, the out-degree of every vertex must me the same, namely
$|S|$, so not every quiver is eligible. A link with cocycle quiver
\[\includegraphics{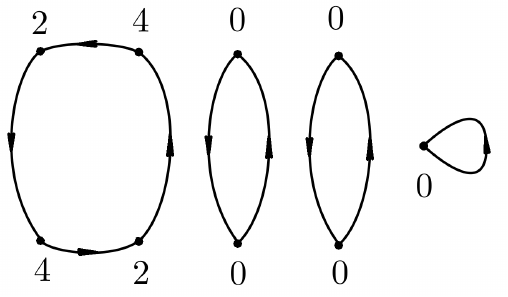}\]
would have the same cocycle quiver polynomial as $L7a3$ but be 
distinguished by the cocycle quiver itself, though we do not know of 
such a link.

To what extent is possible to reverse engineer 
link/quandle/endomorphism/cocycles to fit a particular $R$-colored quiver?
What are necessary and sufficient conditions for a $R$-colored quiver to
be the quandle cocycle quiver of a knot or link, and given such a quiver
how can we construct the set of all links with the given $R$-colored
quiver as $\mathcal{Q}_{X,S}(L)$?

\bibliography{kc-sn2}{}

\begin{thebibliography}{1}

\bibitem{CJKLS}
J.~S. Carter, D.~Jelsovsky, S.~Kamada, L.~Langford, and M.~Saito.
\newblock Quandle cohomology and state-sum invariants of knotted curves and
  surfaces.
\newblock {\em Trans. Amer. Math. Soc.}, 355(10):3947--3989, 2003.

\bibitem{CN1}
K.~Cho and S.~Nelson.
\newblock Quandle coloring quivers.
\newblock {\em J. Knot Theory Ramifications}, 28(1):1950001, 12, 2019.

\bibitem{EN}
M.~Elhamdadi and S.~Nelson.
\newblock {\em Quandles---an introduction to the algebra of knots}, volume~74
  of {\em Student Mathematical Library}.
\newblock American Mathematical Society, Providence, RI, 2015.

\bibitem{J}
D.~Joyce.
\newblock A classifying invariant of knots, the knot quandle.
\newblock {\em J. Pure Appl. Algebra}, 23(1):37--65, 1982.

\bibitem{M}
S.~V. Matveev.
\newblock Distributive groupoids in knot theory.
\newblock {\em Mat. Sb. (N.S.)}, 119(161)(1):78--88, 160, 1982.

\end{thebibliography}
\bibliographystyle{abbrv}

\bigskip

\noindent
\textsc{Department of Mathematics \\
Harvey Mudd College\\
301 Platt Boulevard \\
Claremont, CA 91711
}

\bigskip

\noindent
\textsc{Department of Mathematical Sciences \\
Claremont McKenna College \\
850 Columbia Ave. \\
Claremont, CA 91711}

\end{document}